\documentclass[12pt]{article}
\usepackage{epsfig}
\usepackage{amssymb}
\usepackage{amsfonts}
\usepackage{amsmath}
\usepackage{indentfirst}
\usepackage{epsfig}

\begin{document}

\title{From Fibonacci to Catalan permutations}
\date{}

\author{E. Barcucci, A. Bernini, M. Poneti\\
{\small Dipartimento di Sistemi e Informatica, Universit\`a di Firenze,}\\
{\small Viale G. B. Morgagni 65, 50134 Firenze - Italy}\\
{\tt barcucci, bernini, poneti@dsi.unifi.it}\\
}
\maketitle

\begin{abstract}
It is well known that permutations avoiding any $3$-length pattern
are enumerated by the Catalan numbers. If the three patterns
$123$, $132$ and $213$ are avoided at the same time we obtain a
class of permutations enumerated by the Fibonacci numbers. We
start from these permutations and make one or two forbidden
patterns disappear by suitably ``generalizing'' them. In such a
way we find several classes of permutations enumerated by integer
sequences which lay between the Fibonacci and Catalan numbers. For
each class, we provide the generating function according to the
length of the permutations. Moreover, as a result, we introduce a
sort of ``continuity'' among the number sequences enumerating
these classes of permutations.
\end{abstract}

\section{Introduction}\label{intro}

Fibonacci and Catalan numbers are very well known sequences. They
appear in many combinatorial problems as they enumerate a great
quantity of combinatorial objects. For instance, Fibonacci numbers
are involved in the tiling of a strip, in rabbits' population
growth, in bees' ancestors, \ldots, while Catalan numbers occur in
the enumeration of several kinds of paths, trees, permutations,
polyominoes and other combinatorial structures. Fibonacci numbers
are described by the famous recurrence:
\\
$$ \left\{
\begin{array}{l}
F_{0}=1\\
F_{1}=1\\
F_{n}=F_{n-1}+F_{n-2}\\
\end{array}
\right.
$$
\\
from which the generating function:
\\
$${F(x)} = \frac{1}{1-x-x^2}$$
\\
arises, and the sequence begins with $1,1,2,3,5,8,13,21,34,
\ldots$.
Catalan numbers have been deeply studied, too: they appear in many
relations, also connected to other sequences or by themselves.
They are defined by:
$$
\left\{
\begin{array}{l}
C_{0}=1 \\
C_{1}=1\\
C_n= \sum_{i=0}^{n-1} C_{n-1-i} C_i \\
\end{array}
\right.
$$
\\
The expression
$$
C_n = \frac{1}{n+1}\binom{2n}{n},\quad \mbox{with}\  n\geq0,$$
\\
derived from the generating function
$$
{C(x)} = \frac{1+\sqrt{1-4x}}{2x},$$
\\
is a closed formula for them and the sequence begins with the
numbers $1,1,2,5,14,42,132, \ldots$.\\

Our question is: ``What is there between Fibonacci and Catalan
numbers?'' For instance the following sequences:

\begin{itemize}
    \item $\{c_n\}_{n\geq 0}=\{1,1,2,4,7,13,24,\ldots\}$,
    ($c_0=1,c_1=1,c_2=2,c_n=c_{n-1}+c_{n-2}+c_{n-3}$) Tribonacci
    numbers;
%
    \item $\{t_n\}_{n\geq 0}=\{1,1,2,4,8,16, 2,\ldots,2^{n-1}\}$, ($t_0=1,
    t_n=2^{n-1}$);
    \item $\{p_n\}_{n\geq 0}=\{1,1,2,5,12,29,70,\ldots$\}, ($p_0=1,
    p_1=1, p_2=2, p_n=2p_{n-1}+p_{n-2}$) Pell numbers;
    \item $\{\bar{F}_n\}_{n\geq 0}=\{1,1,2,5,13,34,89,\ldots\}$,
    ($\bar{F}_0=1, \bar{F}_1=1, \bar{F}_n=3\bar{F}_{n-1}+\bar{F}_{n-2}$)
    even index Fibonacci numbers,
\end{itemize}
(for more details see the sequences M1074, M1129, M1413, M1439 in
\cite{Sl}, respectively, where they are defined with different
initial conditions)

\bigskip
\noindent lay between Fibonacci and Catalan numbers (we call the
last sequence \emph{even} index Fibonacci numbers while other
authors call them \emph{odd} index Fibonacci numbers, but this
depends on the initials conditions assumed for the Fibonacci
sequence). We are looking for a unifying combinatorial
interpretation for all these sequences, and others too. To this
aim we will use permutations avoiding forbidden subsequences. Our
results can be extended to paths and trees.

The main idea we are going to base on, has already been used in
\cite{BDPP1}. Here, here we briefly recall that. It is well known
that $|S_n(123,213,312)|=F_n$ and $|S_n(123)|=C_n$, as mentioned
in the abstract. The patterns $213$ and $312$, which are not
present in the second equality, can be seen as particular cases of
more general patterns. More precisely, $213$ can be obtained from
the pattern $r_k=k(k-1)(k-2)\ldots 21(k+1)$ with $k=2$, while
$312$ is the pattern $q_k=1(k+1)k(k-1)\ldots 21$ with $k=2$,
again. When $k$ grows, the patterns $r_k$ and $q_k$ increase their
length, then in the limit ($k$ grows to $\infty$) they can be not
considered in the enumeration of the permutations $\pi$ of
$S_n(123,r_k,q_k)$ since, for each $n\geq 0$, any $\pi$ does not
surely contain a pattern of infinite length. In other words,
starting from the case $k=2$ (involving Fibonacci numbers), for
each $k>2$ we provide a class of pattern avoiding permutations
where the pattern are suitably generalized in order to make them
``disappear'' when $k$ grows, leading to the class $S(123)$
enumerated by the Catalan numbers. We say that there is a sort of
``continuity'' between Fibonacci and Catalan numbers since we
provide a succession of generating functions $\{g_k(x)\}_{k\geq2}$
with $g_2(x)=F(x)$ and whose limit is $C(x)$.

As a matter of fact, in the paper this aim is reached in two
steps: first only the pattern $312$ is generalized so that we
arrive to the class $S(123,213)$ enumerated by $\{2^{n-1}\}$, then
the pattern $213$ is increased in order to obtain the class
$S(123)$. Nevertheless it is possible to make ``disappear'' both
the patterns at the same time obtaining similar results.

\section{Notations and Definitions}\label{notations}

We denote by $S_n$ the set of permutations on $[n]=\{1, 2, \ldots,
n\}$. Let $\pi = \pi_1\pi_2 \cdots \pi_n \in S_n$ and $\Gamma =
\gamma_1\gamma_2 \cdots \gamma_k \in S_k$. We say that $\pi$ does
not contain a subsequence of kind $\Gamma$ (or the pattern
$\Gamma$) if no sequence $j_1 < j_2 < \cdots < j_k$ exists such
that $\pi_{j_i} < \pi_{j_h}$ if and only if $\gamma_i < \gamma_h$.

Let $S_n(\Gamma)$ be the set of permutations not containing a
subsequence of kind $\Gamma$. For instance $7465312 \in S_7(123)$
while $7154326 \notin S_7(123)$ since the subsequence $146$ is of
kind $123$. If $\Gamma_1, \Gamma_2, \ldots, \Gamma_j$ are
permutations we denote by $ S(\Gamma_1, \Gamma_2, \ldots,
\Gamma_j)=S(\Gamma_1)\bigcap S(\Gamma_2)\bigcap \cdots \bigcap
S(\Gamma_j)$ the set of permutations on $[n]$ that do not contain
anyone of the sequences $\Gamma_1, \Gamma_2, \ldots, \Gamma_j$.
For instance $6745231 \in S_7(123, 132, 213)$ while $6475231
\notin S_7(123, 132, 213)$, being the sequence $475$ of kind
$132$.

Permutations avoiding forbidden subsequences have been widely
studied by many authors \cite{BK, BDPP1, BDPP2, Che, Chu, EM, Gi,
Gu, Kra, Kre, Si, St, W1, W2, W3}. A very efficient and natural
method to enumerate classes of permutations was proposed by Chung
et al. \cite{Chu} and Rogers \cite{R}, and, later, by West
\cite{W1}. It consists in generating permutations in $S_n$ from
permutations in $S_{n-1}$ by inserting $n$ in all the positions
such that a forbidden subsequence does not arise (we denote these
positions by a `$\diamond$'). These positions are known as
\emph{active sites}, while a \emph{site} is any position between
two consecutive elements in a permutation or before the first
element or after the last one. If a permutation in
$S_{n-1}(\Gamma_1, \ldots, \Gamma_j)$ contains $k$ active sites,
it generates $k$ permutations in $S_{n}(\Gamma_1, \ldots,
\Gamma_j)$. In the sequel, we denote the $i$-th active site as the
site located before $\pi_i$.

\bigskip
In order to show how we can enumerate classes of permutations by
this method, we consider the class $S_n(123)$. Let $\pi =
\pi_1\pi_2 \cdots \pi_n $ be a permutation in $S_n(123)$ such that
$\pi_1 > \pi_2 > \cdots > \pi_{k-1} < \pi_k$. Then the first $k$
sites are active, since the insertion of $n+1$ in one of these
positions does not create a subsequence of kind $123$. On the
contrary, the sites on the right of $\pi_k$ are not active because
the insertion of $n+1$ produces the subsequence
$\pi_{k-1}\pi_k(n+1)$ which is of kind $123$. Therefore, from the
permutation
\medskip\\
$\diamond\pi_1\diamond\pi_2\diamond\cdots\diamond\pi_{k-1}\diamond\pi_k\pi_{k+1}\cdots\pi_n$
\medskip\\
we obtain the following ones:
\medskip\\
$\diamond(n+1)\pi_1\diamond\pi_2\diamond\cdots\diamond\pi_{k-1}\diamond\pi_k\pi_{k+1}\cdots\pi_n$\\
$\diamond\pi_1\diamond(n+1)\pi_2\cdots\cdots\pi_n$\\
$\diamond\pi_1\diamond\pi_2\diamond(n+1)\pi_3\cdots\pi_n$\\
\vdots\\
$\diamond\pi_1\diamond\pi_2\diamond\cdots\diamond\pi_{k-1}\diamond(n+1)\pi_k\cdots\pi_n$
\medskip\\
which have respectively $(k+1), 2, 3,\ldots, k$ active sites. We
remark that from a permutation $\pi$ having $k$ active sites we
obtain $k$ permutations having $(k+1), 2, 3,\ldots, k$ active
sites, independently from the length of the permutation. Such a
permutation is labelled with $(k)$. We can ``condense'' this
property into a \emph{succession rule} (for more details see
\cite{W2,W3}):
\begin{equation}\label{rscat}
\left\{
\begin{array}{l}
(1) \\
(1)\rightsquigarrow (2)\\
(k)\rightsquigarrow(2)(3)\cdots(k)(k+1)
\end{array}
\right.
\end{equation}
\\
where $(k)\rightsquigarrow(2)(3)\cdots(k)(k+1)$ is the
\emph{production} of a permutation $\pi$ with label $(k)$. The
label $(1)$, said the \emph{axiom} of the succession rule, is the
number of active sites of the empty permutation $\varepsilon$
which is the only permutation with length $n=0$, meaning that
$\varepsilon$ generates the minimal permutation $\pi=1$ with
length $n=1$. In turn, $\pi=\diamond 1 \diamond$ has two active
sites, then it produces two permutations: this fact is described
by the second line of the rule $(1)\rightsquigarrow (2)$ (the
production of the axiom).

The recursive construction of permutation in $S_n(123)$ can also
be represented by a \emph{generating tree}, where each node is a
permutation, the permutations obtained from $\pi$ appear as sons
of $\pi$ and the root is the empty permutation $\varepsilon$ with
length $n=0$. Therefore, on the $n$-th level we have all the
permutations of length $n$ (if we assume the root level is $0$).
The succession rule \ref{rscat} relates the outdegree of each node
in the tree to the outdegree of its sons. Usually, from a
succession rule we can obtain a functional equation or a system of
equations from which one can obtain the generating function
$f(x)=\sum_{n\geq0}a_nx^n$ where $a_n$ is the number of objects on
level $n$. From the above example for $S(123)$, it is possible to
obtain (we omit the calculus) the generating function $C(x)$ for
Catalan numbers. Moreover, $|S_n(123)|=C_n, \ \mbox{for}\ n\geq
0$.

\bigskip
The enumeration of the permutations of $S_n(123, 132, 213)$ is
also briefly illustrated, which is the starting point of our
argument, as recalled in the Introduction. In the permutations of
this class only the first two sites can be active: the insertion
of $n+1$ in another site would produce the subsequence
$\pi_1\pi_2(n+1)$ which is of kind $123$ or $213$. If $\pi_1 <
\pi_2$ then only the first site is active because the insertion of
$n+1$ in the second site would produce the subsequence
$\pi_1\pi_2(n+1)$ which is of kind $132$. Let $\pi = \pi_1\pi_2
\cdots \pi_n $ be a permutation in $S_n(123, 132, 213)$; if $\pi_1
< \pi_2$, from $\diamond\pi_1\pi_2 \cdots \pi_n$ we obtain
$\diamond(n+1)\diamond\pi_1\pi_2 \cdots \pi_n$ which has two
active sites; if $\pi_1 > \pi_2$, from $\diamond\pi_1\diamond\pi_2
\cdots \pi_n$ we obtain $\diamond(n+1)\diamond\pi_1\pi_2 \cdots
\pi_n$  and $\diamond\pi_1(n+1)\pi_2 \cdots \pi_n$ having two and
one active sites, respectively. This construction can be encoded
by the succession rule:
\begin{equation}\label{rsfibo}
\left\{
\begin{array}{l}
(1) \\
(1)\rightsquigarrow(2)\\
(2)\rightsquigarrow(1)(2)
\end{array}
\right.
\end{equation}
The above succession rule is an example of \emph{finite}
succession rule since only a limited number of different labels
appear in it. It is easily seen that it leads to Fibonacci numbers
and $|S_n(123,213,312)|=F_n,\ \mbox{for}\ n\geq 0.$

\bigskip
In the last part of this section, we only note that the
permutations of the class $S(123,213)$, which is the intermediate
step between the above considered classes (see the Introduction),
have exactly two active sites (the first two sites), so that the
corresponding succession rule is
\begin{equation}\label{rs2fin}
\left\{
\begin{array}{l}
(1)\\
(1)\rightsquigarrow(2)\\
(2)\rightsquigarrow(2)(2)
\end{array}
\right.
\end{equation}
It is easy to prove that the related enumerating sequence
$\{t_n\}_{n\geq 0}$ is defined by
$$ \left\{
\begin{array}{l}
t_0=1\\
t_{n}=2^{n-1},\quad n\geq1\\
\end{array}
\right.
$$
and $|S_n(123,213)|=t_n$. The corresponding generating function is
$t(x)=\frac{1-x}{1-2x}$. In the sequel, we refer to this sequence
simply with $\{2^{n-1}\}_{n\geq 0}$.

We conclude by observing that all the considered sequences take
into account the empty permutation which is enumerated by $C_0$,
$F_0$ and $t_0$. Moreover, in each presented succession rule the
axiom refers to it and the production $(1)\rightsquigarrow (2)$
describes its behavior.
%
%

\section{From Fibonacci to $2^{n-1}$}\label{F_n to 2^(n-1) first}

Consider a permutation $\pi\in S_n(123,213,1(k+1)k\ldots 2)$. His
structure is essentially known thanks to \cite{E}, where the
author analyzes the permutations of $S_n(123,132,k(k-1)\ldots
21(k+1))$ which is equivalent to the class we are considering (the
permutations of the former are the reverse complement of the
latter). In the same paper the author shows that those
permutations are enumerated by the sequence of $k$-generalized
Fibonacci numbers, providing also the related generating function.
Here, we give an alternative proof of the same facts by using the
ECO method \cite{BDPP4}. To this aim, we recall the structure of
the permutations referring directly to the class
$S_n(123,213,1(k+1)k\ldots 2)$, nevertheless we omit the easy
proofs that one can recover from \cite{E}.

If $\pi\in S_n(123,213,1(k+1)k\ldots 2)$, then:
\begin{itemize}
    \item either $\pi_1=n$ or $\pi_2=n$;
    \item if $\pi_1=n$, then $\pi=n\tau$, with $\tau\in S_{n-1}(123,213,1(k+1)k\ldots
    2)$;
    \item if $\pi_2=n$, then $\pi_1=n-j$, with $j\in \{1,2,\ldots,
    k-1\}$, and $\pi=(n-j)n(n-1)\ldots (n-j+1)\sigma$, with $\sigma\in
    S_{n-j-2}(123,213,1(k+1)k\ldots 2)$.
\end{itemize}

If $\pi\in S_n(123,213,1(k+1)k\ldots 2)$, denote $\pi^{(i)}$ the
permutations such that $\pi_1=n-i$. The active sites of $\pi$ are
the first two sites: the insertion of $n+1$ in any other site
would create the forbidden pattern $123$ or $213$. More precisely,
the permutations $\pi^{(j)}$ with $j\in\{0,1,2,\ldots,k-2\}$ have
label $(2)$ (the first two sites are active), while $\pi^{(k-1)}$
has label $(1)$ (the first site is active). The son of the
permutation $\pi^{(k-1)}$ is the permutation of
$S_{n+1}(123,213,1(k+1)k\ldots 2)$ obtained from $\pi$ by
inserting $n+1$ in its first active site, which we denote
$\bar{\pi}^{(0)}$. It is easily seen that $\bar{\pi}^{(0)}$ has,
in turn, label $(2)$. The two sons of the permutations with label
$(2)$ are $\bar{\pi}^{(0)}$ and $\bar{\pi}^{(j+1)}$
($\bar{\pi}^{(j+1)}$ is obtained from $\pi$ by inserting $n+1$ in
the second active site). Therefore, all these permutations have,
in turn, label (2) but $\bar{\pi}^{(k-1)}$ which has label (1).
Since all the labels $(2)$ have not the same production, it is
suitable to label each permutation $\pi^{(j)}$
($j\in\{0,1,2,\ldots,k-2\}$) with $(2_j)$ in order to recognize
the permutation $\pi^{(k-2)}$ whose sons have labels $(1)$ and
$(2)$. Then, the above description can be encoded by:
$$
\left\{
\begin{array}{lll}
(1)&& \\
(1)& \rightsquigarrow &(2_0)\\
(2_j)&\rightsquigarrow &(2_0)(2_{j+1}), \quad \quad \mbox{for}\ \ j=0,1,2,\ldots,k-3\\
(2_{k-2})& \rightsquigarrow &(2_0)(1)\\
\end{array}
\right.
$$
\\
We now deduce the generating function $T^k(x,y)$ of the
permutations of $S(123,213,1(k+1)k\ldots 2)$, according to their
length and number of active sites. To this aim we consider the
subsets $T_1$ of the permutations with label $(1)$ and $T_{2_j}$,
with $j=0,1,2,\ldots, k-2$, of the permutations with label
$(2_j)$. It is obvious that these subsets form a partition of
$S(123,213,1(k+1)k\ldots 2)$. Denote with $T_1(x,y)=\sum_{\pi\in
T_1}x^{n(\pi)}y^{f(\pi)}$ the generating function of $T_1$ and
$T_{2_j}(x,y)=\sum_{\pi\in T_{2_j}}x^{n(\pi)}y^{f(\pi)}$ the
generating function of $T_{2_j}$ ($j=0,1,\ldots, k-2$), where
$n(\pi)$ and $f(\pi)$ are the length and the number of active
sites of a permutation $\pi$, respectively. From the above
succession rule the following system is derived:

$$
\left\{
\begin{array}{ccl}
  T_1(x,y) & = & y+xy\sum_{\pi\in T_{2_{k-2}}}x^{n(\pi)} \\
  &&\\
  T_{2_0}(x,y) & = &  xy^2(T_1(x,1)+\sum_{i=0}^{k-2}T_{2_i}(x,1))\\
  &&\\
  T_{2_j}(x,y)& = &xy^2T_{2_{j-1}}(x,1), \quad j=1,2,\ldots, k-2 \qquad .  \\
\end{array}
\right.
$$
\\
Clearly, it is $T^k(x,y)=T_1(x,y)+\sum_{j=0}^{k-2}T_{2_j}(x,y)$
and, if $y=1$, $T^k(x,1)$ is the generating function of the
permutations of $S(123,213,1(k+1)k\ldots 2)$ according to their
length. From the above system (we omit the calculus), it follows:
$$
T^k(x,1)=\frac{1-x}{1-2x+x^{k+1}}\qquad.
$$
Note that if $k$ grows to $\infty$, the generating function $t(x)$
related to the sequence $\{2^{n-1}\}_{n\geq 0}$ (enumerating the
permutations of $S(123,213)$, see Section \ref{notations}) is
obtained. For each $k\geq 2$, we get an expression which is the
generating function of the $k$-generalized Fibonacci numbers. For
$k=1$, the formula leads to $\frac{1}{1-x}$ which is the
generating function of the sequence $\{1\}_{n\geq 0}$ enumerating
the permutations of $S_n(123,213,12)=S_n(12)=n\ (n-1)\ \ldots 2\
1$. For $k=3$ the succession is
$$
\left\{
\begin{array}{lll}
(1)&& \\
(1)& \rightsquigarrow &(2_0)\\
(2_0)&\rightsquigarrow &(2_0)(2_1)\\
(2_1)& \rightsquigarrow &(2_0)(1)\quad ,\\
\end{array}
\right.
$$
which defines the Tribonacci numbers, whose generatring function
is $T^3(x,1)=\frac{1}{1-x-x^2-x^3}$.

\section{From $2^{n-1}$ to Catalan}\label{2^(n-1) to C_n first}

Let $\pi$ be a permutation of $S_n(123,k(k-1)\ldots 21(k+1))$.
Then if $\pi_i=n$ it is $i\in \{1,2,\ldots, k\}$, otherwise if
$\pi_j=n$ with $j\geq k+1$, it should be
$\pi_1>\pi_2>\ldots>\pi_k$ in order to avoid the pattern $123$.
But in this way the entries $\pi_1,\pi_2,\ldots,\pi_k,\pi_j$ are a
pattern $k(k-1)\ldots 21(k+1)$ which is forbidden.

If $\alpha_{\pi}$ denotes the minimum index $j$ such that
$\pi_{j-1}<\pi_j$, we can describe the active sites of $\pi$ by
using $\alpha_{\pi}$.
\begin{enumerate}
    \item If $\alpha_{\pi}=j\leq k$, then the active sites are the
    first $j$ sites of $\pi$. The insertion of $n+1$ in any other
    site would create the pattern $123$. In this case $\pi$ as
    label $(j)$.
    \item If $\alpha_{\pi}>k$, then the active sites of $\pi$ are
    the first $k$ sites since the insertion of $n+1$ in any other
    site would lead to the occurrence of the forbidden patterns
    $k(k-1)\ldots 21(k+1)$ or $123$. In this case $\pi$ has label $(k)$.
\end{enumerate}
In order to describe the labels of the sons of $\pi$, in the
sequel we denote $\bar{\pi}^{(i)}$ the permutation $\bar{\pi}\in
S_{n+1}(123,k(k-1)\ldots 21(k+1))$ obtained from $\pi$ by
inserting $n+1$ in the $i$-th active site of $\pi$.
\begin{enumerate}
    \item If $\pi$ has label $(k)$, it is not difficult to see
    that $\alpha_{\bar{\pi}^{(1)}}=\alpha_{\pi}+1>k$, then
    $\bar{\pi}^{(1)}$ has label $(k)$ again. While, if we consider
    $\bar{\pi}^{(i)}$, with $i=2,3,\ldots, k$, then
    $\alpha_{\bar{\pi}^{(i)}}=i$ and $\bar{\pi}^{(i)}$ has label
    $(i)$. Therefore the production of the label $(k)$ is
    $(k)\rightsquigarrow(2)(3)\ldots(k)(k)$.
    \item If $\pi$ has label $(j)$ with $j\in\{2,3,\ldots,k-1\}$,
    then it is easily seen that $\alpha_{\bar{\pi}^{(1)}}=\alpha_{\pi}+1\leq
    k$ and $\bar{\pi}^{(1)}$ has label $(j+1)$ (note that in this
    case $\alpha_{\pi}=j$). While if we consider
    $\bar{\pi}^{(i)}$, with $i=2,3,\ldots,j$, then
    $\alpha_{\bar{\pi}^{(i)}}=i$ and $\bar{\pi}^{(i)}$ has label
    $(i)$. Therefore the production of $(j)$ is
    $(j)\rightsquigarrow(2)(3)\ldots(j)(j+1)$.
\end{enumerate}
The above construction can be encoded by the succession rule:
$$
\left\{
\begin{array}{lll}
(1)&& \\
(1)& \rightsquigarrow &(2)\\
(j)&\rightsquigarrow &(2)(3)\ldots (j)(j+1), \quad \quad \mbox{for}\ \ j=2,3,\ldots,k-1\\
(k)& \rightsquigarrow &(2)(3)\ldots (k)(k)\qquad ,\\
\end{array}
\right.
$$
where the axiom and its production refer to the empty permutation
\empty generating the permutation $\pi=1$, which, in turn,
produces two sons: $\pi=12$ and $\pi=21$. Using the theory
developed in \cite{DFR}, the production matrix related to the
above succession rule is

$$
P_k=
\left(%
\begin{array}{ccccccc}
  0 & 1 & 0 & \cdots & \cdots & \cdots & \cdots \\
  0 & 1 & 1 & 0 & \cdots & \cdots & \cdots \\
  0 & 1 & 1 & 1 & 0 & \cdots & \cdots \\
  \vdots & \vdots & \vdots & \vdots & \ddots &  &  \\
  \vdots & \vdots & \vdots & \vdots &  & \ddots &  \\
  0 & 1 & 1 & 1 & 1 & \cdots & 1 \\
  0 & 1 & 1 & 1 & 1 & \cdots & 2 \\
\end{array}%
\right)\quad,
$$
with $k$ rows and columns. For each $k\geq 2$, it is easy to see
that the matrix $P_k$ can be obtained from $P_{k-1}$ as follows:

$$
P_k=
\left(%
\begin{array}{cc}
  0 & u^T \\
  0 & P_{k-1}+eu^T \\
\end{array}%
\right)\quad ,
$$
where $u^T$ is the row vector $(1,0,\ldots, 0)$ and $e$ is the
column vector $(1,1,\ldots,1)^T$ (both $k-1$-dimensional). If
$f_{P_k}(x)$ is the generating function according to the length of
the permutations associated to $P_k$, from a result in \cite{DFR}
(more precisely Proposition 3.10), the following functional
equation holds:

$$
f_{P_k}(x)=\frac{1}{1-xf_{P_{k-1}}(x)}\ .
$$

In the limit, we have $f(x)=\frac{1}{1-xf(x)}$ which is the
functional equation verified by the generating function of the
Catalan numbers $C(x)$.
\\
As a particular case, it is possible to check that for $k=3$, the
sequence of the even index Fibonacci numbers is involved. The
obtained succession rule is
$$
\left\{
\begin{array}{lll}
(1)&& \\
(1)& \rightsquigarrow &(2)\\
(2)&\rightsquigarrow &(2)(3)\\
(3)& \rightsquigarrow &(2)(3)(3)\quad ,\\
\end{array}
\right.\
$$
\\
leading to the related generating function
$\bar{F}(x)=\frac{1-2x}{1-3x+x^2}$.

\section{Another street for the same goal}\label{alternative street}

In Section \ref{F_n to 2^(n-1) first}, starting from
$S(123,213,132)$ and using the knowledge that $S(123,213)$ is
enumerated by $\{2^{n-1}\}_{n\geq 0}$, the pattern $132$ has been
generalized in $1(k+1)k\ldots2$, in order to make it ``disappear".
Since the class $S(123,132)$ is enumerated by $\{2^{n-1}\}_{n\geq
0}$, too, one can choose the pattern $213$ instead of $132$ (among
the forbidden patterns of the permutations of $S(123,213,132)$) as
the one to be generalized. Indeed, there is no a particular reason
why we chose the pattern $132$ to make it disappear.

Similarly, starting from $S(123,132)$ and recalling that
$|S_n(p)|=C_n\ \forall\ p\in S_3$, either the pattern $123$ or the
pattern $132$ can be generalized in order to find a class
enumerated by the Catalan numbers.

The difference between a choice with respect to another one lies
in the fact that different ECO construction for the permutations
are expected. Therefore, different succession rules for the same
sequence could be found.

\subsection{From Fibonacci to $2^{n-1}$}\label{F_n to 2^(n-1) second}

Starting from $S(123,213,132)$, here we generalize the pattern
$213$ considering the class $S(123,132,k(k-1)\ldots21(k+1))$, for
$k\geq3$. This class has already been described in \cite{E}, where
the author provides the structure of its permutations. From his
results, it is possible to deduce the following succession rule
(similarly to Section \ref{F_n to 2^(n-1) first}, the details are
omitted), encoding the construction of those permutations:
$$
\left\{
\begin{array}{l}
(1)\\
(1)\rightsquigarrow(2) \\
(h)\rightsquigarrow(1)^{h-1}(h+1)\quad \mbox{for}\ h<k \\
(k)\rightsquigarrow(1)^{k-1}(k)
\end{array}
\right.
$$
In \cite{E} the author shows also that the $k$-generalized
Fibonacci numbers are the enumerating sequence of the permutations
of $S(123,132,k(k-1)\ldots21(k+1))$. This fact can be derived also
by solving the system that can be obtained from the above
succession rule, with a technique similar to that one used in
Section \ref{F_n to 2^(n-1) first} leading to the same generating
function $T^k(x,1)=\frac{1-x}{1-2x+x^{k+1}}$. This agrees with the
fact that in the limit for $k\rightarrow \infty$, the class to be
considered is $S(123,132)$, enumerated by $\{2^{n-1}\}_{n\geq0}$
\cite{Si}. We note that it is possible to describe the
permutations of $S(123,132)$ with the succession rule
$$
\left\{
\begin{array}{l}
(1)\\
(1)\rightsquigarrow(2) \\
(h)\rightsquigarrow(1)^{h-1}(h+1)\quad, \\
\end{array}
\right.
$$
from which one can get that the related generating function is,
again, $t(x)=\frac{1-x}{1-2x}$.

\medskip
The particular case $k=3$ is marked: the obtained succession rule
is
$$
\left\{
\begin{array}{lll}
(1)&& \\
(1)& \rightsquigarrow &(2)\\
(2)&\rightsquigarrow &(1)(3)\\
(3)& \rightsquigarrow &(1)(1)(3)\quad .\\
\end{array}
\right.\
$$
corresponding to the sequence of Tribonacci numbers, as one can
check by deriving the related generating function
$T^3(x,1)=\frac{1}{1-x-x^2-x^3}$.

\subsection{From $2^{n-1}$ to Catalan}\label{2^(n-1) to C_n second}

Starting from $S(123,132)$, the pattern $132$ is generalized in
$(k-1)(k-2)\ldots21(k+1)k$, with $k\geq3$. Moreover, the
construction of the permutations of
$S(123,(k-1)(k-2)\ldots21(k+1)k)$ is described and the
corresponding succession rule is showed. Finally, we prove that
the corresponding generating function is, in the limit for
$k\rightarrow\infty$, the generating function of the Catalan
numbers $C(x)$.

Let $\pi$ be a permutation of $S_n(123,(k-1)(k-2)\ldots21(k+1)k)$.
We denote:
\begin{itemize}
    \item $r=\min\{1,2,\ldots,n\}\ \mbox{such that}\
    \pi_{r-1}<\pi_r$;
    \item $s=\min\{1,2,\ldots,n\}$ and $t=\min\{1,2,\ldots,n\}$ such
    that, fore some indexes $m_1<m_2<\ldots<m_{k-2}<s<t$, it is
    $\pi_{m_1}\pi_{m_2}\ldots\pi_{m_{k-2}}\pi_s\pi_t
    \simeq(k-1)(k-2)\ldots21k$ (the two subsequences are
    order-isomorphic and $\pi_s$ and $\pi_t$ correspond to the $1$ and
    to the $k$ of the pattern $(k-1)(k-2)\ldots21k$);
    \item $\alpha_{\pi}=\min\{r,s\}$;
    \item $\bar{\pi}^{(l)}$ the permutation of
    $S_{n+1}(123,(k-1)(k-2)\ldots21(k+1)k)$ obtained from $\pi$ by
    inserting $n+1$ in the $l$-th site.
\end{itemize}
We prove that $\pi$ has $\alpha_{\pi}$ active sites which are the
first $\alpha_{\pi}$ sites of $\pi$.

It is easily seen that that the insertion of $n+1$ in any site
among the first $\alpha_{\pi}$ sites of $\pi$, does not induce
either the pattern $123$ or the pattern
$(k-1)(k-2)\ldots21(k+1)k$. On the other hand, if
$\alpha_{\pi}=r$, then the insertion of $n+1$ in the $l$-th site,
$l>\alpha_{\pi}$, would create the pattern $123$ in the entries
$\bar{\pi}^{(l)}_{r-1}\bar{\pi}^{(l)}_r\bar{\pi}^{(l)}_l$. While,
if $\alpha_{\pi}=s$, then the insertion of $n+1$ in the $i$-th
site, $\alpha_{\pi}+1\leq i\leq t$, would create the pattern
$(k-1)(k-2)\ldots21(k+1)k$ in the entries
$\bar{\pi}^{(i)}_{m_1}\bar{\pi}^{(i)}_{m_2}\ldots\bar{\pi}^{(i)}_{m_{k-2}}
\bar{\pi}^{(i)}_{\alpha_{\pi}}\bar{\pi}^{(i)}_i\bar{\pi}^{(i)}_{t+1}$
(recall that $\bar{\pi}_i^{(i)}=n+1$ and
$\bar{\pi}_{t+1}^{(i)}=\pi_t$). Finally, if $i\geq t+1$, the
pattern $123$ would appear in the entries
$\bar{\pi}^{(i)}_{\alpha_{\pi}}\bar{\pi}^{(i)}_t\bar{\pi}^{(i)}_i$.

\bigskip\noindent
Denote $(h)$ the label of $\pi$, whit $h=\alpha_{\pi}$. In order
to describe the labels of the sons $\bar{\pi}^{(l)}$,
$l=1,2,\ldots,h$, of $\pi$, we have:
\begin{enumerate}
    \item If $h<k$ (note that on this case $\alpha_{\pi}=r$ or,
    if $\alpha_{\pi}=s$, then $s=k-1$), then
    the permutation $\bar{\pi}^{(1)}=(n+1)\pi_1\pi_2\ldots
    \pi_{\alpha_{\pi}}\ldots\pi_k\ldots\pi_n$,
    so that $\alpha_{{\bar\pi}^{(1)}}=\alpha_{\pi}+1$. Therefore
    $\bar{\pi}^{(1)}$ has label $(h+1)$. While if we consider the
    permutations $\bar{\pi}^{(j)}$, $j=2,3,\ldots,h$, it is
    $\alpha_{\bar{\pi}^{(j)}}=j$ since
    $\bar{\pi}_{j-1}^{(j)}<\bar{\pi}_j^{(j)}(=n+1)$. So
    $\bar{\pi}^{(j)}$ has label $(j)$ and we conclude that the
    production of $(h)$ is
    $(h)\rightsquigarrow(2)(3)\ldots(h)(h+1)$.
    \item If $h\geq k$, then $\bar{\pi}^{(1)}=(n+1)\pi_1\pi_2\ldots
    \pi_k\ldots\pi_{\alpha_{\pi}}\ldots\pi_n$, so that
    $\alpha_{{\bar\pi}^{(1)}}=\alpha_{\pi}+1$. Therefore $\bar{\pi}^{(1)}$
    has label $(h+1)$. Note that in both cases $\alpha_{\pi}=r$ or
    $\alpha_{\pi}=s$ it is
    $\pi_1>\pi_2>\ldots>\pi_{\alpha_{\pi}-1}$. Then, if we
    consider the permutations $\bar{\pi}^{(j)}$,
    $j=k,k+1,\ldots,\alpha_{\pi}$, we obtain $\alpha_{\bar{\pi}^{(j)}}=k-1$, regardless of
    $j$, since
    $\bar{\pi}^{(j)}_1\bar{\pi}^{(j)}_2\ldots\bar{\pi}^{(j)}_{k-1}\bar{\pi}^{(j)}_j
    \simeq(k-1)(k-2)\ldots 1k$. Then $\bar{\pi}^{(j)}$ has label
    $(k-1)$, for $j=k,k+1,\ldots,\alpha_{\pi}$. For the remaining
    sons $\bar{\pi}^{(j)}$, $j=2,3,\ldots,k-1$, it is easily seen that
    $\bar{\pi}^{(j)}_{j-1}<\bar{\pi}^{(j)}_j(=n+1)$. So,
    $\bar{\pi}^{(j)}$ has label $(j)$. We conclude that, in this
    second case, the production of $(h)$ is $(h)\rightsquigarrow(2)(3)
    \ldots(k-2)(k-1)^{h-k+2}(h+1)$.
\end{enumerate}
The above description of the generation of the permutations of
$S(123,(k-1)(k-2)\ldots21(k+1)k)$ can be then encoded in the
following succession rule $\Omega_k$:
$$\Omega_k=
\left\{
\begin{array}{l}
(1) \\
(1)\rightsquigarrow(2)\\
(h)\rightsquigarrow(2)\cdots(h)(h+1)\quad \mbox{for}\ h<k\\
(h)\rightsquigarrow(2)\cdots(k-2)(k-1)^{h-k+2}(h+1)\quad \mbox{for
}\ h\geq k\qquad .
\end{array}
\right.
$$
For $k=2$, the class $S(123,132)$ is obtained, whose corresponding
succession rule has been considered in Section \ref{F_n to 2^(n-1)
second}. Note that it does not correspond with the one obtained
from the above one poising $k=2$.

For $k=3$ (the class is $S(123,2143)$) we get the succession rule:
$$
\left\{
\begin{array}{l}
(1)\\
(1)\rightsquigarrow(2) \\
(h)\rightsquigarrow(2)^{h-1}(h+1)\quad ,
\end{array}
\right.
$$
leading to the even index Fibonacci numbers. Note that it is
different from the succession rule corresponding to the same
numbers of Section \ref{2^(n-1) to C_n first}. Its associated
production matrix \cite{DFR} is: \scriptsize
$$
M_3=\left(%
\begin{array}{cccccc}
  0 & 1 & 0 & 0 & 0 & \cdots \\
  0 & 1 & 1 & 0 & 0 & \cdots \\
  0 & 2 & 0 & 1 & 0 & \cdots \\
  0 & 3 & 0 & 0 & 1 & \cdots \\
  \vdots & \vdots & \vdots & \vdots & \vdots & \ddots \\
\end{array}%
\right)\ .
$$
\normalsize For each $k\geq4$, it is easy to check that the
production matrix related to $\Omega_k$ satisfies
$$
M_k=
\left(%
\begin{array}{cc}
  0 & u^T \\
  0 & M_{k-1}+eu^T \\
\end{array}%
\right)\quad ,
$$
where $u^T=(1,0,0,\ldots)$ and $e=(1,1,1,\ldots)^T$. Then, if
$g_{M_k}(x)$ is the corresponding generating function, we deduce
\cite{DFR}:
$$
g_{M_k}(x)=\frac{1}{1-xg_{M_{k-1}(x)}}\ .
$$
If $g(x)$ denotes the limit of $g_{M_k}(x)$, the functional
equation $g(x)=\frac{1}{1-xg(x)}$ is obtained, which is verified
by the generating function $C(x)$ of the Catalan numbers.

\section{From Fibonacci to Catalan directly}\label{F_n to C_n directly}

This section summarizes the results found when the two patterns
$132$ and $213$ are generalized at the same time, considering the
class $S(123,(k-1)(k-2)\ldots21(k+1)k),k(k-1)\ldots21(k+1))$ in
order to obtain the class $S(123)$, when $k$ grows to $\infty$.
Most of the proofs are omitted but they can easily recovered by
the reader. At the first step, for $k=3$, we find the succession
rule:
$$
\left\{
\begin{array}{l}
(1)\\
(1)\rightsquigarrow(2) \\
(2)\rightsquigarrow(2)(3)\\
(3)\rightsquigarrow(2)(2)(3)
\end{array}
\right.
$$
\\
corresponding to $S_n(123, 2143, 3214)$. This class is enumerated
by Pell numbers which we define with the recurrence:
$$
\left\{
\begin{array}{l}
p_{0}=1\\
p_{1}=1\\
p_2=2\\
p_{n}=2p_{n-1}+p_{n-2},\quad n\geq3
\end{array}
\right.
$$
\\
Note that the initial conditions are different from the usual ones
(which are $p_0=0$ and $p_1=1$) in order to consider the empty
permutation $\varepsilon$, for $n=0$.

For a general $k$ we have the class $S_n(123, (k-1)\cdots1(k+1)k,
k(k-1)\cdots1(k+1))$. We briefly describe the construction of the
permutations of the class (the details are omitted). Let $\pi$ be
a permutation of the class. It is easily seen that if $\pi_l=n$,
then $l\leq k$. Therefore, if $(h)$ denotes the label of $\pi$, it
is $h\in\{1,2,\ldots,k\}$. Now, if $h<k$, then $\bar{\pi}^{(1)}$
has label $(h+1)$ and $\bar{\pi}^{(j)}$, $j=2,3,\ldots,h$, has
label $(j)$. While, if $h=k$, then $\bar{\pi}^{(1)}$ has label
$(k)$, $\bar{\pi}^{(j)}$, $j=2,3,\ldots,k-1$, has label $(j)$ and
$\bar{\pi}^{(k)}$ has label $(k-1)$, again. The construction can
be encoded in the succession rule:
$$
\left\{
\begin{array}{l}
(1)\\
(1)\rightsquigarrow(2) \\
(h)\rightsquigarrow(2)(3)\cdots(h-1)(h)(h+1) \mbox{    }h<k\\
(k)\rightsquigarrow(2)(3)\cdots(k-1)(k-1)(k)\quad .
\end{array}
\right.
$$
For each $k$, considering the associated production matrices
\cite{DFR} and the corresponding generating functions, it possible
to prove that, in the limit, the generating function of the
Catalan numbers is obtained.

\subsection{A continuity between Pell numbers and even index Fibonacci numbers}

We conclude by showing that it is possible to find a
``continuity'' between Pell and even index Fibonacci numbers. We
start from the class $S_n(123, 2143, 3214)$ (obtained by posing
$k=3$ in the preceding succession rule) enumerated by Pell
numbers, then we generalize the pattern $2143$, so obtaining the
classes $S(123, 3214, 21(k+1)k (k-1)\ldots43)$.

Let $\pi\in S_n(123, 3214, 21(k+1)k (k-1)\ldots43)$. Then, if
$\pi_l=n$, it is $l\leq3$ in order to avoid the patterns $123$ and
$3214$. Therefore, $\pi$ has at most $3$ active sites (the first
three sites of $\pi$). We denote $r_{\pi}$ the number of entries
of $\pi$ with index $j\geq 3$ such that $\pi_j>\pi_1$ (note that
if $\pi_1>\pi_2$, then $r_{\pi}=0$). It is:
\begin{itemize}
    \item $\pi_{j_1}>\pi_{j_2}>\ldots>\pi_{j_{r_{\pi}}}$ (the pattern $123$ is
    forbidden);
    \item $r_{\pi}\leq (k-2)$ (the pattern $21(k+1)k\ldots43$ is
    forbidden);
    \item the elements $\pi_{j_i}$ are adjacent in $\pi$ in order
    to avoid $123$ or $21(k+1)k\ldots43$.
\end{itemize}

If $\pi$ starts with an ascent (i.e. $\pi_1<\pi_2$), then only the
first two sites are active, since the insertion of $n+1$ in any
other site would create the pattern $123$: the permutation $\pi$
has label $(2)$.

If $\pi$ starts with a descent (i. e. $\pi_1>\pi_2$), then the
number of its active sites depends on $r_{\pi}$:
\begin{enumerate}
    \item If $r_{\pi}=h<k-2$, then $\pi$ has three active sites.
    Let $(3_h)$ be its label. The permutation $\bar{\pi}^{(1)}$
    (obtained by $\pi$ by inserting $n+1$ in the first site)
    starts with a descent and $r_{\bar{\pi}^{(1)}}=0$ (since
    $\bar{\pi}^(1)_1=n+1$); therefore, $\bar{\pi}^{(1)}$ has label
    $(3_0)$. The son $\bar{\pi}^{(2)}$ starts with an ascent and
    its label is $(2)$. The last son $\bar{\pi}^{(3)}$ starts with
    a descent and $r_{\bar{\pi}^{(3)}}=h+1$, so its label is
    $(3_{h+1})$. The production of $(3_h)$ is
    $(3_h)\rightsquigarrow(2)(3_0)(3_{h+1})$.
    \item If $r_{\pi}=k-2$, then $\pi$ has two active sites, since
    the insertion in the third site would create the pattern
    $21(k+1)k\ldots43$, while the insertion in any other site
    surely creates the pattern $123$. Its son $\bar{\pi}^{(1)}$
    has label $(3_0)$ since it starts with a descent and
    $r_{\pi^{(1)}}=0$. While the other son $\bar{\pi}^{(2)}$
    starts with an ascent and has label $(2)$. Therefore, the
    production of label $(2)$ is $(2)\rightsquigarrow(2)(3_0)$.
\end{enumerate}
The following succession rule:
$$
\left\{
\begin{array}{l}
(1)\\
(1)\rightsquigarrow(2)\\
(2)\rightsquigarrow(2)(3_0)\\
(3_j)\rightsquigarrow(2)(3_0)(3_{j+1}),\quad \mbox{for}\ \ j=0,1,2,\ldots,k-3\\
(3_{k-3})\rightsquigarrow(2)(2)(3_0)\\
\end{array}
\right.
$$
summarizes the construction of the class
$S(123,3214,21(k+1)k\ldots43)$. Solving the system one can deduce
from the above rule, the generating function
$\bar{F}_k(x)=\frac{1-2x+x^k}{1-3x+x^2+x^k}$ is obtained, which in
the limit is the generating function of the even index Fibonacci
numbers $\bar{F}(x)$.

\bigskip
Starting from the class $S(123,2143,3214,)$, one can generalize
the pattern $3214$ instead of $2143$. The class we get is
$S(123,2143,k(k-1)\ldots32(k+1)1)$ and the succession rule
describing its construction is (the easy proof is omitted):
$$
\left\{
\begin{array}{l}
(1) \\
(1)\rightsquigarrow(2)\\
(h)\rightsquigarrow(2)^{h-1}(h+1)\quad \mbox{for }\ h<k\\
(k)\rightsquigarrow(2)^{k-1}(k)\quad.
\end{array}
\right.
$$
Once again, one can prove that the corresponding generating
function is $\bar{F}_k(x)$, leading, in the limit, to
$\bar{F}(x)$.



\section{Remarks}
\begin {figure}[h]
\begin {center}
\includegraphics [scale=0.4] {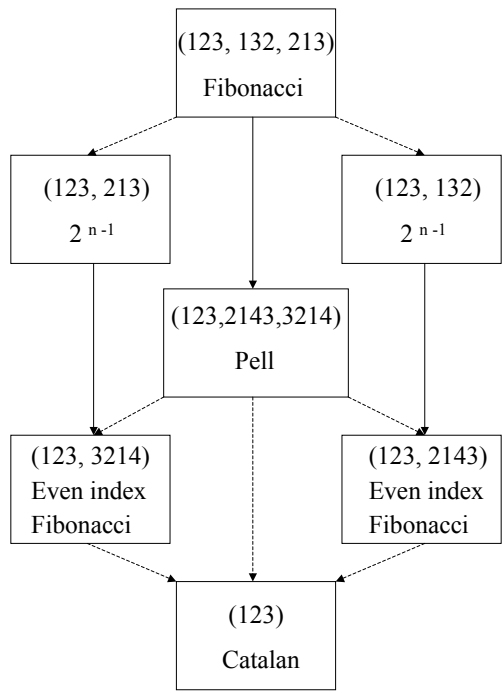}
\end {center}
\caption {The graph of permutations.}\label{fig1}
\end {figure}

In order to summarize the several ``continuities'' we have
proposed in the paper, we condense our results in Figure
\ref{fig1} where a straight line represents a direct step and a
dashed line represents a family of permutations obtained by
generalizing one or two patterns.

\bigskip
The results we found for permutations can be easily extended to
Dyck paths and planar trees by means of ECO method
\cite{BDPP3,BDPP4}. We can find classes of paths and trees
described by the finite succession rules we introduced by imposing
some conditions on the height of paths and the level of their
valleys and on the outdegree and level of nodes in the trees.

\bigskip
Figure \ref{fig1} allows to see the different three ways we have
followed to describe a discrete ``continuity" between Fibonacci
and Catalan numbers: the generalization of a single pattern (the
rightmost and the leftmost path from the top to the bottom in the
figure) and the generalization of a pair of patterns (central path
in the figure). In particular, following the rightmost and the
leftmost path in the graph, the intermediate level of the
permutations enumerated by $\{2^{n-1}\}_{n\geq0}$ is encountered.
For each $k$, our approach produces two different class of
permutations enumerated by the same sequence, indeed the two
corresponding generating functions are the same for each $k$. We
note that, in this way, we can provide two different succession
rules encoding the same sequence. An instance can be seen by
looking at the succession rules the reader can find at the end of
the Sections \ref{F_n to 2^(n-1) first} and \ref{F_n to 2^(n-1)
second}.

The same happens with the succession rule at the end of Section
\ref{2^(n-1) to C_n second} and the succession rule of the
particular case ($k=3$) of Section \ref{2^(n-1) to C_n second},
which encode the sequence of the even index Fibonacci numbers.
Really, we did not prove that this is the case for each $k$
related to the classes of permutations used to describe the
discrete continuity between $\{2^{n-1}\}_{n\geq0}$ and Catalan
numbers, since we did not get the explicit formulas of the
generating functions.


\end{document}